# SWENDSEN-WANG DYNAMICS ON $\mathbb{Z}^d$ FOR DISORDERED NON FERROMAGNETIC SYSTEMS.

EMILIO DE SANTIS

Abstract. We study the Swendsen-Wang dynamics for disordered non ferromagnetic Ising models on cubic subsets of the hypercubic lattice $\mathbb{Z}^d$ and we show that for all small values of the temperature parameter $\beta^{-1}$ the dynamics has a slow relaxation to equilibrium (it is *torpid*). Looking into this dynamics from the point of view of the Markov chains theory we can prove that the spectral radius goes to one when the size of the system goes to infinity. This means that, if we want to use the Swendsen-Wang dynamics for a computer simulation, we have a slow convergence to the stationary measure in low temperature. Also it is a good example of a non-local dynamics that relaxes slowly to the equilibrium measure.

## 1. Introduction

The study of the Swendsen-Wang (SW) dynamics has received a considerable attention in recent years; see the original paper [26] and also [15, 19, 20, 7], where some rigorous mathematical results on the convergence speed to the equilibrium measure have been obtained.

The Swendsen-Wang algorithm was proposed for the simulation on a computer of the Ising or Potts model (not necessarily ferromagnetic) to obtain a more rapid simulation with respect to the standard Monte Carlo methods. In this paper we are precisely interested in the behavior of the algorithm for non ferromagnetic Ising models defined on cubic subsets of $\mathbb{Z}^d$. We do not deal with general Potts model, but we conjecture that many of the results proved here are true, in some sense, also in that case.

The SW dynamics is a reversible Markov chain having the Gibbs Ising measure as the stationary one. A relevant feature of the SW algorithm is that it is non-local, in the sense that many spins can change value simultaneously. This feature allows us to hope for a fast convergence to equilibrium. In fact for ferromagnetic systems with an external field on $\mathbb{Z}^d$ it is proved (see [19, 20]) that this method of simulation has a fast approach to the equilibrium measure; it is effective also for small values of the temperature parameter where standard Monte Carlo methods are typically slow. On the other hand, it has been recently proved that the speed of convergence to the equilibrium measure of the SW algorithm for the mean field is not *always* rapid [15]. Recently it was proved that the mixing time of the SW dynamics for the Potts model in a box of $\mathbb{Z}^d$ of side length $L$, at the critical temperature, with a large enough number of states, goes to infinity as $L \to \infty$; see [2] where some explicit inequalities are established for the mixing time of the SW algorithm and for the mixing time of the more standard Glauber dynamics.

We will prove that the SW dynamics for *disordered non ferromagnetic Ising systems*, on $\mathbb{Z}^d$, spin-glass systems, has a slow approach to the equilibrium measure for large values of the parameter $\beta$ (the inverse of the temperature). More precisely we prove that the spectral gap of the associated Markov chain decreases to zero when the side length of the box goes to infinity if $\beta$ is large enough.

Some arguments and results of the paper can be framed in the dynamical Griffiths singularity. The dynamical Griffiths singularity for the ordinary Glauber dynamics has been studied by various authors (see [1, 9, 5, 6, 12, 16]). Notice also that our results have some similarities with [2], in the sense that we show another model (but with a large region of the parameter temperature in which we can apply the proofs) where the SW algorithm has a slow convergence to equilibrium for the physically relevant case of the lattice $\mathbb{Z}^d$.

We notice that in [18] it is proved that the Swendsen-Wang dynamics for ferromagnetic Ising model in the low temperature case is non-ergodic and hence there is a slow convergence to the equilibrium in finite volume. In any case there are no results on the spectral gap and we conjecture that the spectral gap remains bounded away from zero as the side length goes to infinity for the SW dynamics applied to the ferromagnetic Ising model. Thus we believe that it is not easy to replace our construction to prove that the spectral gap shrinks to zero. Furthermore, we conjecture that the the behavior is not governed by the occurrence of large ferromagnetic blocks, almost isolated, as it is in some other cases (see for example [5, 6, 8]). Also for these reasons we believe that this kind of Griffiths singularity is a new phenomenon not studied yet.

Our results are true in any dimension for all the Ising models where the interactions are random variables and 0 belongs to the interior of the support. This hypothesis on the support of the measure can seem unnatural but there are a lot of probability measures that respect those conditions; for example the uniform measure on an interval that includes the origin or the Gaussian measures, so the results are highly applicable.

We stress that the SW dynamics is related to the Fortuin-Kasteleyn (FK) representation that has an important role to study the phase transition. In fact for ferromagnetic models there is a 1-1 relationship between bond percolation and spin phase transition (see [14, 21]) but for frustrated systems there is a weaker relationship that allows us only to establish inequalities on the phase transition for the spin-glass (see [4, 8, 13, 14, 21]). In any case it is not known for the general spin-glass if there is a phase transition or not. Thus our results are not proved using arguments on the equilibrium measures or the phase transition of these Ising models, but we use directly the structure of the SW dynamics.

## 2. Swendsen-Wang dynamics and main definitions

In this section we will recall the definition of the Swendsen-Wang dynamics and the main definitions that will be used in the paper. Our goal is to prove that the spectral radius associated with the Markov chain representing the SW dynamics goes to one when the size of the system grows to infinity invading all the lattice $\mathbb{Z}^d$. In order to show the essential line of the proof, we will not find explicit bounds on the spectral radius but we will only prove that it goes to one. In any case all the estimations could be done explicitly. So, if one is interested in usable bounds for the spectral radius, it is possible to perform explicit calculations following our proofs.

Let us consider the lattice $\mathbb{Z}^d$ where the edges are all the couples of vertices $x$, $y$ having $||x-y|| = 1$ and $||\cdot||$ is the Euclidean distance.

*Some sets.* Given a set of edges $B \subset E$ we define the *vertex set of* $B$, $V(B)$, as the set of all vertices $v \in V$ such that there exists $w \in V$ with $\{v, w\} \in B$. Given a set of vertices $A$ we define $E(A)$ as the set of edges $e \in E$ such that $V(\{e\}) \in A$; we remark that given a set of vertices $V_1 \subset V$ we have $V(E(V_1)) \subset V_1$ and given a set of edges $E_1 \subset E$ the relation $E(V(E_1)) = E_1$ is true. In the following if we do a set operation using a vertex set $A$ and an edge set $B$ in the same time, then we will see these two sets as vertex sets; for example if $A \subset V$ and $B \subset E$ then $A \cap B$ corresponds to $A \cap V(B)$.

*Boundaries.* For a set of vertices $V \subset \mathbb{Z}^d$ let us denote by $\partial V$, the *boundary of* $V$, as the set of vertices belonging to $V$ such that for every vertex $u \in V$ there is a vertex $w \in V^c$ such that $\{u, w\} \in E$ is an edge.

For a set of vertices $V$ we define also the *edge boundary* $\partial_b V$ as follows $e = \{u, v\} \in \partial_b V$ *if and only if* $e \cap V \neq \emptyset$ and $e \cap V^c \neq \emptyset$.

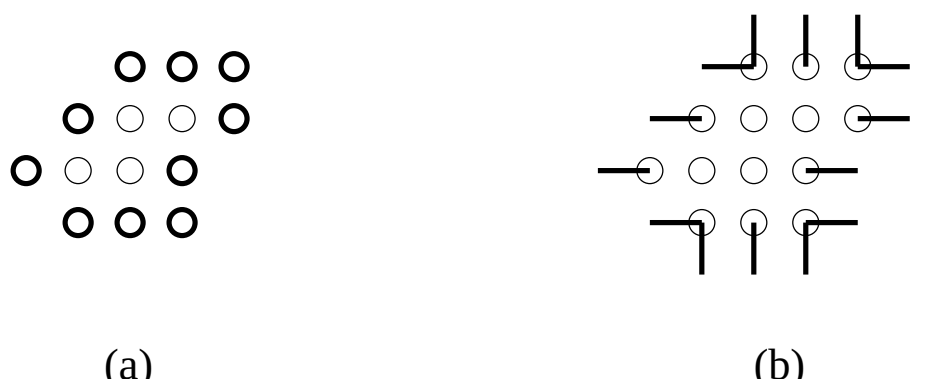

FIGURE 1. The boundary of a set of vertices in (a) and its edge boundary in (b)

We will deal with interactions $\mathcal{J} = \{J_b\}_{b \in E(\mathbb{Z}^d)}$ that are *i.i.d.* random variables with a distribution $Q_J$ having a neighborhood of zero inside its support (supp $Q_J$).

Our proof works for all $\beta \in (\beta_h, \infty)$ where $\beta_h$ depends on the distribution $Q_J$ and on the dimension $d$ of the lattice. We could also generalize the proof to some other kinds of distributions of the interactions, for example the uniform measure on $(-b, -a) \cup (a, b)$ for suitable $a, b > 0$, but only for some limited intervals of the parameter $\beta$ and not for all $\beta \in (\beta_h, \infty)$.

To fix the ideas and the notation we can think that $Q_J$ is a continuous distribution whose support includes $[-1, 1]$. If this is not the case we can re-scale the support of $Q_J$ such that $[-1, 1] \subset$ supp $Q_J$; this will correspond to change the value of the parameter $\beta$.

We will consider an *(hyper) cube of vertices* $B_{4l}(v)$ with *length of the side* $4l$ and *center* $v$. Let us fix an event in the realization of interaction on the cube $B_{4l}(v)$ having positive probability.

For the edges $e \in \partial_b B_{4l}$ we require that $\sum_{e\in\partial_b B_{4l}} |J_e| < 1/l^2$; for every edge $e \in E(B_{4l}) \setminus \partial_b B_{2l}$ we set $J_e \in (1 - 1/l^d, 1]$; finally for every $e \in \partial_b B_{2l}$ we define the interactions through these rules (only for $l$ large enough):

1. We choose $d$ faces $S_1, S_2, \ldots, S_d$ among the $2d$ faces of the cube $B_{2l}$ in such a way that every face $S_i$ is not the image of another of these faces using $R$, the *central reflection of the cube.* The faces $S_i$ are sets of vertices. The central reflection of the cube $R$ can be used as a function on the vertices of the cube into itself or as a function on the edges of the cube into itself (it is an isomorphism of the graph). So for every edge $e = \{v_1, v_2\}$ of the cube also $Re := \{Rv_1, Rv_2\}$ is an edge.
2. For every edge $e \in \partial_b B_{2l}$ having a vertex $v = (v_1, v_2, \ldots, v_d)$ belonging to a face $\{S_i\}_{i=1,\ldots,d}$ of the cube we give the following rule: if $|\sum_{i=1}^d v_i|$ is an even number we impose that $J_e \in (a_d - 1/l^d, a_d)$ otherwise $J_b \in (-sa_d, -sa_d + 1/l^d)$ where $s$ is in $(0, 1)$ (in the following we will specify it) and $a_d$ is an opportune positive constant.
3. For every edge $e \in \partial_b B_{2l}$ that we have not already defined we use this rule: if $J_{Re} \in (a_d - 1/l^d, a_d)$ then $J_e \in (-sa_d, -sa_d + 1/l^d)$ and if $J_{Re} \in (-sa_d, -sa_d + 1/l^d)$ then $J_e \in (a_d - 1/l^d, a_d)$.

So all the interactions inside the cube $B_{4l}$ are defined.

We have fixed the interactions on a finite region so that there is a positive probability to find a cube $B_{4l}$ of this kind. When are fixed all the interactions on the graph $\mathbb{Z}^d$ one can find a cube $B_{4l}$ of this kind at a finite distance from the origin $Q_J$-almost surely (by ergodicity with respect to the translation). We will indicate with $\Lambda_L$ the box $[-L, L]^d \cap \mathbb{Z}^d$. Now, using standard arguments, it is easy to see that there exists a subsequence of boxes $\Lambda_{L(l)}$ with $L(l) \in \mathbb{N}$ strictly increasing with $l$ such that only for a finite number of $l$ the box $\Lambda_{L(l)}$ does not contain a cube $B_{4l}$ with the desired configuration of interactions ($Q_J$-almost surely).

*Configuration Space.* Let us define the space of the spin configurations $\Omega = \{-1,1\}^{\mathbb{Z}^d}$; let us write $\Omega_\Lambda = \{-1,1\}^\Lambda$ where $\Lambda$ is a finite subset of $\mathbb{Z}^d$. We define the space $H = \{0,1\}^{E(\mathbb{Z}^d)}$ (configurations of edges) and $H_\Lambda = \{0,1\}^{E(\Lambda)}$ (configurations of edges on $\Lambda$).

*Ising model.* Let us define the Hamiltonian of the Ising model on a finite subset $\Lambda \subset \mathbb{Z}^d$ as:

$$H_{\Lambda,\mathcal{J}}(\sigma) = -\sum_{i,k:||i-k||=1} J_{i,k}\sigma_i\sigma_k. \tag{1}$$

Here $i,k$ are vertices in $\Lambda \cap \mathbb{Z}^d$ and $\{i,k\} \in E(\Lambda)$, $\sigma = \{\sigma_i\}_{i\in\Lambda\cap\mathbb{Z}^d}$ is a configuration belonging to $\Omega_\Lambda$.

The Ising model (or Ising measure) on a finite subset $\Lambda \subset \mathbb{Z}^d$ is defined as

$$\pi_{\Lambda,\mathcal{J}\beta}(\sigma) = \frac{e^{-\beta H_{\Lambda,\mathcal{J}}(\sigma)}}{Z_{\Lambda,\mathcal{J}\beta}} \tag{2}$$

where $Z_{\Lambda,\mathcal{J}\beta}$ is a normalizing factor (*partition function*). Sometimes we will omit the dependence from the box $\Lambda$ and the parameters $\mathcal{J}$.

Note that the interactions $\mathcal{J} = \{J_e\}_{e\in E(\Lambda)}$ are random variables; in this case the system is called *disordered* or *quenched* since there will be a random measure or a stochastic process parameterized by some random variables.

*Swendsen-Wang dynamics.* The Swendsen-Wang dynamics (or process) is a Markov chain that updates alternatively the configuration of edges $H_\Lambda$ and the configuration of spins $\Omega_\Lambda$. The SW dynamics has as space of states $\Omega_\Lambda \times H_\Lambda$ with a finite $\Lambda$ (we can consider a cube in $d$ dimensions as $\Lambda$) and it is an algorithm generally used to simulate the Ising model on a finite space.

The Markov transition probability $P_\beta(\omega,\eta;\tilde{\omega},\tilde{\eta})$, with $\omega,\tilde{\omega} \in \Omega_\Lambda$ and $\eta,\tilde{\eta} \in H_\Lambda$, that realizes the SW dynamics is not homogeneous, indeed in every odd time $t \in \mathbb{N}$ the algorithm updates the configuration of the edges (in $H_\Lambda$) and in every even time $t \in \mathbb{N}$ SW dynamics updates the configuration of the spins (in $\Omega_\Lambda$). We will also introduce the Markov transition probability $M_\beta^{SW}(\omega,\tilde{\omega})$, with $\omega,\tilde{\omega} \in \Omega_\Lambda$, that is the Markov transition probability on the even times restricted to the space of the states $\Omega_\Lambda$ (so $M_\beta^{SW} = P_\beta^2|_{\Omega_\Lambda}$). For a configuration $(\sigma,\eta) \in \Omega_\Lambda \times H_\Lambda$ at time $t = 2n$, with $n \in \mathbb{N}$, the algorithms update the configuration with these two steps:

SW1) *This step updates the configuration of edges.* Let $\tilde{A} \subset E(\Lambda)$ be the subset of edges $e = \{x, y\} \in E(\Lambda)$ such that $J_e\sigma_x\sigma_y \geq 0$. For every edge $e \in \tilde{A}$, independently of the others, with probability $p_e = 1 - e^{-\beta|J_e|}$ in the updated configuration $\eta'_e = 1$ and with probability $1 - p_e$ in the updated configuration $\eta'_e = 0$ (we remark that if $J_e = 0$ then $p_e = 0$).

SW2) *This step updates the configuration of vertices.* Let us consider the e-clusters, maximal connected component, of the graph where two vertices, $u$ and $w$, are connected if there exists a path $\gamma = (u, e_1, \ldots, e_n, w)$ such that $\eta_{e_i} = 1$ for every $i = 1, \ldots, n$. For every e-cluster $C$ we will choose a vertex $v$ and we will decide with probability $1/2$ that $\sigma_v$ is equal to 1 or $-1$. For all the other vertices of the e-cluster the values of the spin are assigned following the rule $J_e\sigma_x\sigma_y \geq 0$ where $e = \{x, y\}$. This procedure define in every vertex the configuration $\Omega_\Lambda$.

It is immediate to see that the Markov chain $M^{SW}_\beta$ is irreducible and aperiodic (since every element of the transition matrix $M^{SW}_\beta$ is larger than zero); $\pi_\beta$, the Gibbs measure of the Ising model, is the invariant measure of the Markov chain $M^{SW}_\beta$ (see Edwards and Sokal [26, 11]). We will indicate with $\nu_\beta : \Omega_\Lambda \times H_\Lambda \to [0, 1]$ the invariant measure of the Markov chain $P_\beta$. The measure $\nu_\beta$ is the FK-measure related to the Ising model (see [21]).

It is easy to understand that also the transition probability $P_\beta$ is irreducible and aperiodic on the states $\Omega_\Lambda \times H_\Lambda$. It is believed that the SW algorithm is faster than the standard Monte Carlo algorithm because this dynamics is not local so it is more easy to leave a metastable state [26].

*Reversible property.* The SW algorithm is a *reversible* Markov chain; it was proved in the original paper [26]. The reversibility or the *detailed balance* was directly used in [26] to show that the equilibrium measure of this algorithm is the Ising measure. So a consequence of reversibility and regularity (i.e. the Markov chain is irreducible and aperiodic) is that all the eigenvalues of the transition matrix are real and inside the interval $(-1, 1]$. These eigenvalues have multiplicity one with probability one if the interactions $\{J_e\}_{e\in E}$ are absolutely continuous random variables. We order these eigenvalues

$$1 = \lambda_0 > \lambda_1 > \cdots > \lambda_n > -1 \tag{3}$$

where $n$ is the number of configurations on a finite number of vertices $V$ that is equal to $2^{|V|}$.

*Exponential Time and the Spectral Radius.* For the exposition of the *exponential time* we follow [24]. We will give the definitions for a system with a finite number of states. Let $P$ be a regular transition matrix associated to a Markov chain and let $\pi$ be its unique invariant measure. Let $f : \Omega_\Lambda \to \mathbb{R}$ be a real function that is square-integrable with respect to the stationary measure $\pi$ , i.e. $f \in L^2(\pi)$.

Let us write

$$\pi(f) = \sum_{\omega \in \Omega} f(\omega)\pi(\omega) \tag{4}$$

for the mean value of $f$ with respect to the measure $\pi$ and

$$C_{ff}(t) = \sum_{\omega,\tilde{\omega} \in \Omega} f(\omega)f(\tilde{\omega})\pi(\omega)P^t(\omega,\tilde{\omega}) - [\pi(f)]^2 = \tag{5}$$

$$= Var_\pi(f) - \frac{1}{2}\sum_{\omega,\tilde{\omega} \in \Omega} (f(\omega) - f(\tilde{\omega}))^2\pi(\omega)P^t(\omega,\tilde{\omega}) = Var_\pi(f) - \mathcal{E}(f,f) \tag{6}$$

for the *autocorrelation function.* We are denoting with $Var_\pi(f)$ the variance of the random variable $f$ and with $\mathcal{E}_t(f,f)$ the *Dirichlet form* of the transition matrix $P^t$; in the following the Dirichlet form will be used to bound the spectral radius of the transition matrix. We also notice that $C_{ff}(0) = Var_\pi(f)$.

Let us define the *normalized autocorrelation function* as

$$\rho_{ff}(t) = \frac{C_{ff}(t)}{C_{ff}(0)}. \tag{7}$$

For finite regular Markov chains $\rho_{ff}(t)$ decays exponentially for large $t$; we define, for $f \in L^2(\pi)$,the *exponential f-autocorrelation time*

$$\tau_{exp,f} = \limsup_{t\to\infty} \frac{t}{-\log|\rho_{ff}(t)|} \tag{8}$$

and the *exponential autocorrelation time*

$$\tau_{exp} = \sup_{f \in L^2(\pi)} \tau_{exp,f}. \tag{9}$$

Thus, $\tau_{exp}$ is the typical time that the Markov chain employs to relax towards the equilibrium measure.

Equivalently one can see $P$ as an operator on the Hilbert space $L^2(\pi)$. Let define $R$ the *spectral radius* of the operator $P$ as

$$R = \sup_{i>0} |\lambda_i| \tag{10}$$

If there is only the eigenvalue $\lambda_0 = 1$ on the unit circle then $R$ for a finite space of states is strictly less than 1. The relation

$$R = e^{-1/\tau_{exp}} \tag{11}$$

is in [ST89, Propositions 2.3-2.5] (see also [24]). So if $R$ is close to one the Markov chain relax slowly to the equilibrium measure.

Our Markov chain is dependent from the size of the box $\Lambda_L = [-L, L]^d \cap \mathbb{Z}^d$, so the quantities $\tau_{exp,f} = \tau_{exp,f,L}$, $\tau_{exp} = \tau_{exp,L}$, and $R = R_L$ are function of $L$. We will choose some particular functions $f_L$ to show, using the variational characterization, that $R_L$ goes to one almost surely (see Theorem 4.8) for $L$ that goes to infinity (and hence $\tau_{exp,L} \to \infty$).

*Magnetization and related events.* Let us define the *magnetization* of a configuration $\sigma \in \Omega_\Lambda$ on a finite region $A \subset \Lambda$ of vertices as:

$$M(A, \sigma) = \frac{1}{|A|} \sum_{v \in A} \sigma_v, \tag{12}$$

so $M(A, \sigma)$ is a rational number belonging to the interval $[-1, 1]$.

Let us define the *accordance on the set of edges* $B \subset E(\Lambda)$ as:

$$K(B, \sigma) = \frac{1}{|B|} \sum_{e=\{v,w\} \in B} \sigma_v \sigma_w. \tag{13}$$

Finally let's define the events:

$$S_l^+(\delta) = \{\sigma : K(\partial_b B_{2l}, \sigma) \geq 1 - \delta\}, \tag{14}$$

$$S_l^-(\delta) = \{\sigma : K(\partial_b B_{2l}, \sigma) \leq -1 + \delta\}, \tag{15}$$

where $\delta$ is a constant in $(0,1)$.

It will not be difficult to see that for all $\delta \in (0,1)$ and $\beta$ large enough the probability of the event $S_l^+(\delta)$ tends to one when $l$ goes to infinity. In Section 4 we will give some detailed estimations for the probability of the events $S_l^+(\delta)$ and $S_l^-(\delta)$ that will be used to prove that the spectral radius goes to one when $\Lambda \to \mathbb{Z}^d$.

## 3. Satisfied edges and Peierls surfaces

In this section we will give other definitions.

*Clusters.* For every $v \in \mathbb{Z}^d$ let us denote by $C_v = C_v(\sigma)$ the *cluster at the vertex* $v$. $C_v$ is defined as the maximal connected set of $\mathbb{Z}^d$ having constant sign. In this case two vertices $u$ and $w$ are connected if there is a path $\gamma = (u, e_1, v_1, \ldots, v_{n-1}, e_n, w)$ such that $\sigma_u = \sigma_{v_i}$ and $\sigma_u = \sigma_w$ for every $i = 1, \ldots, n-1$. In an obvious way we can extend this definition of cluster to a generic graph $\mathcal{G} = (V, E)$.

Given a set of vertices $V$ (resp. a set of edges $E$) we will call $C$ *a cluster on* $V$ (resp. on $E$) the cluster on the subgraph $\mathcal{H} = (V, E(V))$ (resp. $\mathcal{H} = (V(E), E)$).

*Peierls surface.* We call a set of edges $S$ a *Peierls surface* if there exists a configuration $\sigma$ such that $S = \partial_b C(\sigma)$, where $C$ is a cluster. It is possible to associate a *plaquette* to every edge $e$ of the lattice, where this plaquette is the $(d-1)$-dimensional unit cube orthogonal to $e$ intersecting the edge in its center. It is a standard fact that the plaquettes associated to a Peierls surface form a closed (without boundary) $(d-1)$-manifold.

*Satisfied, unsatisfied and neutral edges.* For a fixed configuration $\sigma$ we call *satisfied* an edge $e = \{u, v\}$ such that $J_{uv}\sigma_u\sigma_v > 0$. It is called *unsatisfied* if $J_{uv}\sigma_u\sigma_v < 0$ and *neutral* if $J_{uv}\sigma_u\sigma_v = 0$. It is easy to see that given a set of vertices $V$ if we consider a spin-flip on $V$ (i.e. the sign of the spin of all the vertices in $V$ changes) then all the edges $e \in \partial_b V$ that were satisfied (resp.unsatisfied) become unsatisfied (resp.satisfied); all the other edges that are not in $\partial_b V$ or that are neutral remain of the same kind. Note that if $J$ has an absolutely continuous distribution then the probability that there is a neutral edge is zero.

## 4. Spectral radius for the SW dynamics

In the first two lemmas of this section we present some geometrical results. Lemma 4.1 is an easy calculation that shows an inclusion between two events. The following Lemmas 4.2, 4.3 are a consequence of some results on the minimal surfaces; we write the proofs following the reference [3]. In Lemmas 4.4, 4.6 and 4.7 we give the essential bounds on the probability that the Markov chain does a transition from any state in $S_l^-(\delta)$ to any state in $S_l^+(\delta)$ and a bound on the probability of the event $\Omega \setminus (S_l^-(\delta) \cup S_l^+(\delta))$. In Theorem 4.8 we combine the previous bounds and, using the variational characterization of the second eigenvalue of the transition matrix, we will show that the spectral radius goes to one so the exponential autocorrelation time goes to infinity. In all the bounds of the lemmas we do not use that the Markov chain is reversible but in Theorem 4.8 it is essential for the use of the variational characterization.

**Lemma 4.1.** *Let $\delta \in (0, 1/2)$. If $\sigma \in \Omega \setminus (S_l^+(\delta) \cup S_l^-(\delta))$ and $l$ is large enough then at least one of the following inequalities is true:*

a) $|M(\partial B_{2l}, \sigma)| < 1 - \delta/4$;

b) $|M(\partial B_{2(l+1)}, \sigma)| < 1 - \delta/4$.

**Proof**. The proof of this Lemma is an algebraic calculation. Let us suppose that $M(\partial B_{2l}, \sigma) \geq 0$. If $M(\partial B_{2l}, \sigma) < 0$ we can do a global spin flip of the configuration leaving unchanged the absolute values of the magnetization and the value of $K$. So we can consider only the case $M(\partial B_{2l}, \sigma) \geq 0$. Now there are two different cases to analyze:

a) $M(\partial B_{2(l+1)}, \sigma) \geq 0$;

b) $M(\partial B_{2(l+1)}, \sigma) < 0$.

First we prove the lemma in case a). As a consequence of our hypothesis $K(\partial_b B_{2l}, \sigma) < 1 - \delta$ so there are at least $\delta |\partial_b B_{2l}|/2$ edges $e = \{v, w\}$ on which the product of $\sigma_v \sigma_w$ is equal to $-1$. So at least $\delta |\partial B_{2l}|/4$ are the negative spins that are either on $\partial B_{2l}$ or on $\partial B_{2(l+1)}$. Therefore we have that either

$M(\partial B_{2l},\sigma) < 1-\delta/2$ or

$$M(\partial B_{2(l+1)},\sigma) < \frac{1}{|\partial B_{2(l+1)}|}\left[|\partial B_{2(l+1)}| - 2\frac{\delta}{4}|\partial B_{2l}|\right].$$

It is clear that $\lim_{l\to\infty}|\partial B_{2l}|/|\partial B_{2(l+1)}| = 1$ so for $l$ large enough $|\partial B_{2l}|/|\partial B_{2(l+1)}| > 1/2$ and we obtain the inequality of the lemma. Case b) follows by an analogous calculation just using the inequality $K(\partial_b B_{2l},\sigma) > -1+\delta$. □

The following two lemmas are connected to problems of minimal surfaces for a given boundary.

**Lemma 4.2.** *Let $\delta \in (0,1/2)$. If for a configuration $\sigma \in \Omega$ the inequality $|M(\partial B_{2l},\sigma)| < 1-\delta$ is verified then there exists a Peierls surface $\gamma$ inside box $B_{2l}$ where all the edges in $\gamma\setminus\partial_b B_{2l}$ are unsatisfied and $|\gamma\setminus\partial_b B_{2l}| > \delta\rho_d|\partial_b B_{2l}|$, where $\rho_d$ is a constant in $(0,1]$ depending only on the dimension $d$.*

We will use also the following analogous lemma that we will not prove because the proof is basically the same of Lemma 4.2.

**Lemma 4.3.** *Let $\delta \in (0,1/2)$. If for a configuration $\sigma \in \Omega$ the inequality $|M(\partial B_{2(l+1)},\sigma)| < 1-\delta$ is verified then there exists a Peierls surface $\gamma$ inside $B_{4l}\setminus B_{2l}$ where all the edges in $\gamma\setminus(\partial_b B_{2l}\cup\partial_b B_{4l})$ are unsatisfied and $|\gamma\setminus(\partial_b B_{2l}\cup\partial_b B_{4l})| > \delta\rho_d|\partial_b B_{2l}|$, where $\rho_d$ is a constant in $(0,1]$ depending only on the dimension $d$.*

Before the proof of Lemma 4.2 we give these two definitions

$$A_l^+(\sigma) = \{v\in\partial B_{2l} : \sigma_v = 1\}, \tag{16}$$

$$A_l^-(\sigma) = \{v\in\partial B_{2l} : \sigma_v = -1\} = \partial B_{2l}\setminus A_l^+(\sigma). \tag{17}$$

We notice that, as an easy consequence of hypothesis of Lemma 4.2, there exists a constant $\rho = \min\{\frac{\delta}{2},\frac{1}{10}\} \in (0,\frac{1}{10}]$ such that

$$\rho < \frac{|A_l^+(\sigma)|}{|\partial B_{2l}|} < 1-\rho \tag{18}$$

**Proof of Lemma 4.2**. To fix ideas we impose $0 \leq M(\partial B_{2l}, \sigma) < 1-\delta$ but the other case $-1+\delta < M(\partial B_{2l}, \sigma) < 0$ is analogous. Let us consider the region $B_{2l}$ and the family of clusters $\mathcal{C} = \{\tilde{C}_i\}_{i=1,\dots,M}$ on this region. The family of clusters $\mathcal{C}$ form a partition of the vertices in $V(B_{2l})$; we fix an order on the clusters belonging to $\mathcal{C}$. We will do a spin-flip of all the vertices belonging to the first cluster $C \in \mathcal{C}$ which do not touch the boundary $\partial B_{2l}$; sometimes in the following, doing this spin-flip of the vertices of a cluster $C$, we will say that we are eliminating the cluster $C$. We will fix a new order on this new family of clusters on $B_{2l}$ and we will continue to eliminate the first cluster that do not touch the boundary $\partial B_{2l}$. In the end of this procedure we will obtain a configuration $\sigma^{(a)}$ in which all the clusters on $B_{2l}$ touch the boundary $\partial B_{2l}$. We will denote with $\mathcal{C}^{(a)} = \{C_i(\sigma^{(a)})\}_{i=1,\dots,M^{(a)}}$ the new family of clusters on $B_{2l}$. It is easy to verify that the boundary of the clusters can only decrease with this procedure, so if we can find the Peierls surface $\gamma$ (verifying Lemma 4.2) in this new configuration $\sigma^{(a)}$ this is allowed also in the old configuration $\sigma$. So we will find a Peierls surface $\gamma$ such that:

$$\gamma \subset \cup_i \partial_b C_i(\sigma^{(a)}) \subset \cup_i \partial_b C_i(\sigma). \tag{19}$$

Also it is easy to verify that $A_l^+(\sigma) = A_l^+(\sigma^{(a)})$, and so $M(\partial B_{2l}, \sigma^{(a)}) = M(\partial B_{2l}, \sigma)$. Now we will analyze two sub-cases

i) $\partial C_O \cap \partial B_l = \emptyset$,

ii) $\partial C_O \cap \partial B_l \neq \emptyset$,

where $O$ is the center of the cube $B_l$. In the first case we consider $\gamma = \partial_b C_O$ and we do the projection from the point $O$ of the surface $\gamma$ on $\partial B_{2l}$, we obtain the inequality of the Lemma with standard arguments (see [3]).

Indeed in this first case we have:

$$\frac{|\gamma \cap \partial_b B_{2l}|}{|\partial_b B_{2l}|} \leq 1-\rho \tag{20}$$

in fact either $\partial C_O \subset A^+(\sigma^{(a)})$ or $\partial C_O \subset A^-(\sigma^{(a)})$.

Before we analyze case $ii)$ let us define the following partial order for the clusters on the cube $B_{2l}$: the cluster in the center of the cube $C_O$ is the smallest one ($C_O \preceq C_x$ for every $x \in B_{2l}$) and two clusters $C_x$, $C_y$ are ordered $C_x \preceq C_y$ if for every path $\alpha$ starting in the origin and ending in a vertex of $C_y$ the path $\alpha$ has at least a vertex belonging to the cluster $C_x$.

If there is a cluster $\tilde{C}$ such that

$$\frac{\rho}{10} \leq \frac{|\partial\tilde{C} \cap \partial B_{2l}|}{|\partial B_{2l}|} \leq 1 - \frac{\rho}{10} \tag{21}$$

we define $\gamma = \partial\tilde{C}$, otherwise we start the following procedure. Let us choose a maximal cluster with respect to the partial order and consider the spin-flip of this cluster. We notice that with this procedure we can only decrease the boundary of the clusters inside the square $B_{2l}$ so after every such step up is more difficult to find a surface $\gamma$ verifying the Lemma; so if we can prove that there is such a kind of surface then there is the same surface also in the original configuration. Every time that we do a spin-flip of a cluster we will verify if there exist a cluster respecting the inequality (21). When there is a cluster $\tilde{C}$ verifying the inequality (21) we will define $\gamma = \partial\tilde{C}$. Surely there exists a time in which this procedure ends because the relation (21) can be not true only if for all the clusters

$$\frac{|\partial\tilde{C} \cap \partial B_{2l}|}{|\partial B_{2l}|} < \frac{\rho}{10}. \tag{22}$$

So in the end of this procedure there is a unique cluster $\tilde{C}$ that verify the inequalities

$$\frac{\rho}{10} < \frac{|\partial\tilde{C} \cap \partial B_{2l}|}{|\partial B_{2l}|} < \frac{\rho}{5}. \tag{23}$$

Let us consider the intersection of the cluster $\tilde{C}$ with $\partial B_m$ where $m \in \{l, l+1, \ldots, 2l\}$. Let us define

$$\rho_m := \frac{|\tilde{C} \cap \partial B_m|}{|\partial B_m|} \qquad m \in \{l, l+1, \ldots, 2l\}; \tag{24}$$

let's also write

$$\nu := \inf_{m,n} \frac{\rho_m}{\rho_n} \qquad m,\ n \in \{l, l+1, \ldots, 2l\}. \tag{25}$$

We remark that by definition $\nu \leq 1$. Again we will analyze two sub-cases

a) $0 \leq \nu \leq \frac{9}{10}$;

b) $\frac{9}{10} < \nu \leq 1$.

The sub-case $a)$ can be proved with the projection of the surface as in the case $i)$ of the previous page (just changing the constant $\rho_d$).

We prove the sub-case $b)$; on every level $m \in \{l, \dots, 2l\}$ it is easy to see that $\frac{\rho}{10} < \rho_m < 1 - \frac{\rho}{10}$. Now using the isoperimetric inequality on a compact manifold, see [3] pag. 209, for every level $m \in \{l, \dots, 2l\}$ we obtain the inequality

$$|\partial\tilde{C} \cap \partial B_m| > K_d \left[\frac{\rho}{10}|\partial B_m|\right]^{\frac{d-2}{d-1}}, \tag{26}$$

where $K_d$ is a positive constant depending on the manifold (in this case only on the dimension of the space because the manifolds are cubes). For dimension $d = 2$ the RHS of the inequality (26) becomes independent of $\partial B_m$ and it is immediate to understand that

$$|\partial\tilde{C} \cap \partial B_m| > 1. \tag{27}$$

Therefore we have a lower bound for the measure of the surface $S = (\partial\tilde{C} \setminus E(B_l)) \subset \gamma$, $S$ is the surface $\partial\tilde{C}$ out of the square $B_l$ ($S$ is $(d-2)$ dimensional manifold but it is also a set of edges); so we can write

$$K_d l \inf_{m \in \{l,\dots,2l\}} \left[\frac{\rho}{10}|\partial B_m|\right]^{\frac{d-2}{d-1}} > \frac{1}{10} K_d \rho d^{\frac{d-2}{d-1}} l^{d-1}. \tag{28}$$

In the last inequality we use that $\rho \in (0,1]$ so $\rho^{\frac{d-2}{d-1}} > \rho$ for every $d \geq 2$ and $|\partial B_l| = 2dl^{d-1}$. This concludes the proof for the last sub-case, so the lemma is proved for all $d \geq 2$.

□

We fix the value of the constant $a_d$ as:

$$a_d = \rho_d \delta/2 \tag{29}$$

but every positive smaller value is fine; we remind that $a_d$ is approximatively the absolute value of some interactions on $\partial B_{2l}$, so in general $a_d$ is a small positive number.

Given a set of edges $E$ let us define $H_E$ as the *contribution to the Hamiltonian on the edges* $E$; precisely we have:

$$H_E(\sigma) = -\sum_{\{x,y\}\in E} J_{x,y}\sigma_x\sigma_y. \tag{30}$$

For the next lemma see also [10, 9].

**Lemma 4.4.** *Let $\delta \in (0, 1/2)$. There exists $\beta_1(\delta)$ such that for all $\beta > \beta_1(\delta)$ the following upper bound is satisfied*

$$\pi_\beta(\{|M(\partial B_{2l},\sigma)| < 1-\delta\}) + \pi_\beta(\{|M(\partial B_{2(l+1)},\sigma)| < 1-\delta\}) < e^{-\tilde{K}_1\beta l^{d-1}} \tag{31}$$

*where* $\tilde{K}_1 = d2^{d-1}\delta^2\rho_d^2$.

**Proof**. If the absolute value of the magnetization is less than $(1-\delta)$ on the set $\partial B_{2l}$ or on the set $\partial B_{2(l+1)}$ then, using Lemma 4.2 and Lemma 4.3, there exists a Peierls surface $\gamma$ where at least a percentage $\delta\rho_d$ of edges are unsatisfied; moreover all these edges have interactions that are close to 1. The satisfied edges are all on $\gamma \cap (\partial_b B_{2l} \cup \partial_b B_{4l})$ and the absolute value of the interactions of these edges is small (the values are close to $a_d$ or $sa_d$ if $e \in \partial_b B_{2l}$ and close to zero if $e \in \partial_b B_{4l}$). For a fixed Peierls surface $\gamma$ if a configuration $\sigma$ respects the previous conditions about the satisfied and unsatisfied edges then we write $\sigma \sim \gamma$. To bound the probability we consider a spin flip on the interior of $\gamma$. This changes the Hamiltonian only on the edges of $\gamma$ (see fig. 4) leaving the same values on all the other edges; this gives us an easy bound on the probability of all the configurations $\sigma \sim \gamma$ (see [10, 9]). On the edges of

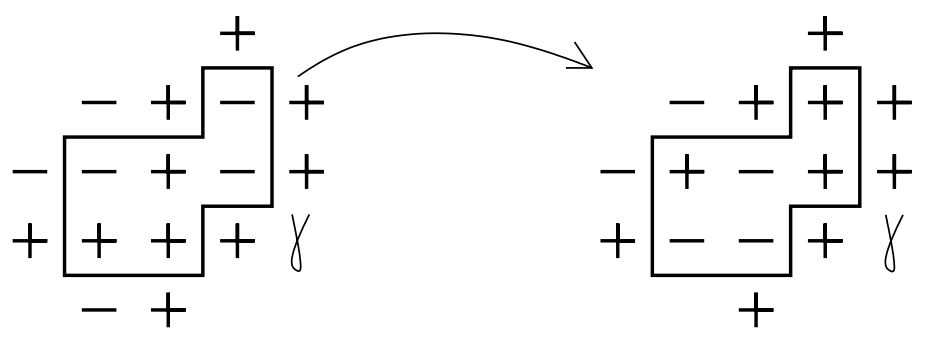

FIGURE 2. The spin-flip of the region inside the circuit $\gamma$

$\gamma = \gamma(\sigma)$ -where $\gamma$ is the Peierls surface of Lemma 4.2 and 4.3- the value of the Hamiltonian has these

bounds:

$$\delta\rho_d|\gamma| - a_d|\gamma| - 1 \leq H_\gamma(\sigma) \leq \delta\rho|\gamma| + a_d|\gamma| + 1, \tag{32}$$

for all $\sigma$ verifing $|M(\partial B_{2l}, \sigma)| < 1 - \delta$ or $|M(\partial B_{2(l+1)}, \sigma)| < 1 - \delta$. The inequalities are clear because all the edges in $\gamma \setminus (\partial_b B_{2l} \cup \partial_b B_{4l})$ -they are at least $\delta\rho|\gamma|$- are unsatisfied; in the first inequality we are considering that all the other edges in $\gamma \cap (\partial_b B_{2l} \cup \partial_b B_{4l})$ are satisfied so they can give a contribute to the Hamiltonian that can be only larger than $-a_d|\gamma| - 1$; vice versa in the second inequality we are considering that also the edges in $\gamma \cap (\partial_b B_{2l} \cup \partial_b B_{4l})$ are unsatisfied. Analogously, taking the spin-flip of the configuration on the interior of $\gamma$, we obtain

$$-\delta\rho_d|\gamma| - a_d|\gamma| - 1 \leq H_\gamma(f_\gamma\sigma) \leq -\delta\rho|\gamma| + a_d|\gamma| + 1, \tag{33}$$

So

$$H(\sigma) - H(f_\gamma\sigma) = H_\gamma(\sigma) - H_\gamma(f_\gamma\sigma) \geq 2\delta\rho_d|\gamma| - 2a_d|\gamma| - 2. \tag{34}$$

Now, using relation (29), one has

$$H(\sigma) - H(f_\gamma\sigma) \geq \delta\rho_d|\gamma| - 2. \tag{35}$$

Let us fix a Peierls surface $\gamma$ such that $|\gamma \setminus (\partial_b B_{2l} \cup \partial_b B_{4l})| \geq \delta\rho_d|\gamma|$. For all the $\sigma$ having all the edges on $\gamma \setminus (\partial_b B_{2l} \cup \partial_b B_{4l})$ unsatisfied we find this bound (we are using inequality (35))

$$\frac{\pi_\beta(\sigma)}{\pi_\beta(f_\gamma\sigma)} = \frac{e^{-\beta H(\sigma)}}{e^{-\beta H(f_\gamma\sigma)}} \leq e^{-\beta(\delta\rho_d|\gamma|-2)}. \tag{36}$$

For a fixed $\gamma$ we perform the sum over all the $\sigma \sim \gamma$ (i.e. all the configurations $\sigma$ having $\gamma \setminus (\partial_b B_{2l} \cup \partial_b B_{4l})$ as a subset of boundary of clusters):

$$\sum_{\sigma\sim\gamma} \pi_\beta(\sigma) \leq e^{-\beta(\delta\rho_d|\gamma|-2)} \sum_{\sigma\sim\gamma} \pi_\beta(f_\gamma\sigma) \leq e^{-\beta(\delta\rho_d|\gamma|-2)}, \tag{37}$$

since the function $f_\gamma : \Omega \to \Omega$ is injective, and so the sum $\sum_{\sigma\sim\gamma} \pi_\beta(f_\gamma\sigma)$ is a probability (therefore a number belonging to $[0, 1]$). Now, taking all the possible $\gamma$ which fulfill the conditions of Lemma 4.2-4.3, and using inequality (37), we obtain:

$$\pi_\beta(\{|M(\partial B_{2l},\sigma)|<1-\delta\})+\pi_\beta(\{|M(\partial B_{2(l+1)},\sigma)|<1-\delta\})$$

$$\leq \sum_{\gamma:|\gamma|\geq \delta\rho_d|\partial_b B_{2l}|} e^{-\beta(\delta\rho_d|\gamma|-2)} \leq \sum_{n\geq \delta\rho_d|\partial_b B_{2l}|} K_0^n e^{-\beta(\delta\rho_d n-2)}, \tag{38}$$

where $K_0^n$ is an uniform upper bound for the number of the Peierls surface of size $n$ in any dimension $d$ (see [17, 22]). Now, if $\beta$ is large enough, we can bound the (38) with the expression

$$\exp\left(-\beta\delta^2\rho_d^2|\partial B_{2l}|/2\right)=\exp\left(-\beta\delta^2\rho_d^2 2^{d-1}d(l^{d-1}+\rho(l)\right),$$

where $\rho(l)$ is a $o(l^{d-1})$. So for $\beta$ large enough we can choose $\tilde{K}_1 = d2^{d-1}\delta^2\rho_d^2$. □

**Remark 4.5.** *If the Markov chain starts with the stationary measure then also the event* $(S_l^+(4\delta)\cup S_l^-(4\delta))^c$ *has probability smaller than* $e^{-K_1\beta l^{d-1}}$ *at every time* $t\in\mathbb{N}$ *(we are also using Lemma 4.1 to connect the magnetization with the events* $S_l^+(4\delta)$ *and* $S_l^-(4\delta)$*). So we have:*

$$\pi_\beta\left((S_l^+(\delta)\cup S_l^-(\delta))^c\right) < e^{-K_1\beta l^{d-1}} \tag{39}$$

*where* $K_1=\tilde{K}_1/4 = d2^{d-3}\delta^2\rho_d^2$

Let us define the following function $\Phi:\Omega\to\Omega$ inside a cube $B_{4l}$ with interactions verifying the conditions 1, 2 and 3 of the second section:

$$\Phi(\sigma)(v)=\begin{cases} -\sigma(Rv) & \text{if } v\in B_{2l}\\ \sigma(Rv) & \text{if } v\in (B_{4l}\setminus B_{2l})\\ \sigma(v) & \text{if } v\notin B_{4l}\end{cases} \tag{40}$$

and $R$ is the reflection with respect to the center of the cube $B_{4l}$. It is easy to see that $\Phi^2=\mathbf{1}$ (the identity function). By construction we have that, for every $\delta\in(0,1)$, if $\sigma\in S^+(\delta)$ then $\Phi(\sigma)\in S^-(\delta)$ and if $\sigma\in S^-(\delta)$ then $\Phi(\sigma)\in S^+(\delta)$. The function $\Phi$ has these important properties:

$$M(\partial B_{2l},\Phi(\sigma)) = -M(\partial B_{2l},\sigma) \text{ and } M(\partial B_{2(l+1)},\Phi(\sigma))=M(\partial B_{2(l+1)},\sigma). \tag{41}$$

Also every spin inside the cube $B_{4l}$ sees the same next neighbored except the spins on the boundaries $\partial B_{2l}$, $\partial B_{2(l+1)}$ and $\partial B_{4l}$. We notice that an edge $e \in \partial_b B_{2l}$ is satisfied for the configuration $\sigma$ if and only if $Re$ is satisfied for the configuration $\Phi(\sigma)$; this observation will be essential in the proof of the next lemma.

**Lemma 4.6.** *Let's $\delta \in (0, 1/2)$. Then there exists $s \in (0, 1-2\delta)$ and $\beta_2(\delta)$ such that for all $\beta > \beta_2(\delta)$:*

$$a) \quad \lim_{l\to\infty} \pi_\beta(S_l^+(\delta)) = 1, \tag{42}$$

$$b) \ \limsup_{l\to\infty} \frac{\pi_\beta(S_l^+(\delta))}{\pi_\beta(S_l^-(\delta))e^{K_2\beta l^{d-1}}} = 0, \tag{43}$$

*where the positive constant $K_2$ is equal to $d2^d(1-s-2\delta)\delta\rho_d$.*

**Proof**. By Lemma 4.4 and Remark 4.5 we know that the probability of $S_l^+(\delta) \cup S_l^-(\delta)$ goes to one when $l$ goes to infinity. Now we will prove that the measure will be concentrated on the set $S_l^+(\delta)$ for $l$ tending to infinity.

Using the function $\Phi : \Omega \to \Omega$ we want to compare the probability of any configuration $\sigma \in S_l^+(\delta)$ with the probability of the configuration $\Phi(\sigma) \in S_l^-(\delta)$. For every realization of the interactions that respects our constraints (see Section 2 points 1-3) we have an uncertainty about the Gibbs factor of order $e^{\pm\beta C}$, where $C$ is a constant independent of $l$, because the contribution to the Hamiltonian inside the cube $B_{4l}$ depends on the fact that the interactions are values belonging to a small interval of length $1/l^d$ (the inverse of the volume so they are inversely proportional to the number of interactions inside the cube $B_{4l}$). We do not control the value of the interactions on the boundary of the faces (they have dimension $(d-2)$ so the number of interactions are $O(l^{d-2})$). Therefore, including the two kinds of *mistakes*, we can calculate the Hamiltonian with an uncertainty of order $l^{d-2}$. Let us consider the configuration $\Phi(\sigma)$ where $\sigma$ is any configuration in $S_l^+(\delta)$. As noted, just before this lemma, we have that every satisfied (resp. unsatisfied) edge $e \in \partial_b B_{2l}$ for the configuration $\sigma$ corresponds to a satisfied (resp. unsatisfied) edge $Re \in \partial_b B_{2l}$ for the configuration $\Phi(\sigma)$. Moreover, if $e \in \partial_b B_{2l}$ has $J_e > 0$ (resp. $J_e < 0$), then $Re$ has $J_{Re} < 0$ (resp. $J_{Re} > 0$), this is a consequence of the definition of the interactions inside the cube $B_{4l}$. Now we give some rigorous upper and lower bounds for the expression $\pi_\beta(\Phi(\sigma))/\pi_\beta(\sigma)$ when $\sigma \in S_l^+(\delta)$; we will use some other bounds for the the Hamiltonian on the edge set $\partial_b B_{2l}$. Doing algebraic calculations we find:

$$\left(\frac{1}{2} - \frac{s}{2} - \delta\right) a_d|\partial_b B_{2l}| \leq H_{\partial_b B_{2l}}(\sigma) \leq \left(\frac{1}{2} - \frac{s}{2} + s\delta\right) a_d|\partial_b B_{2l}| \qquad \forall \sigma \in S_l^-(\delta, t); \tag{44}$$

where the lower (upper) bound is obtained with a configuration $\sigma$ having $K(\partial_b B_{2l}, \sigma) = 1 - \delta$ with all the positive (negative) edges $e = \{i, j\} = \partial_b B_{2l}$ having $\sigma_i \sigma_j = 1$. We also find:

$$\left(-\frac{1}{2} + \frac{s}{2} - s\delta\right) a_d|\partial_b B_{2l}| \leq H_{\partial_b B_{2l}}(\sigma) \leq \left(-\frac{1}{2} + \frac{s}{2} + \delta\right) a_d|\partial_b B_{2l}| \qquad \forall \sigma \in S_l^+(\delta). \tag{45}$$

So, if the parameter $\delta$ is small, all the allowed values of $H_{\partial_b B_{2l}}$ on $S_l^-(\delta)$ (resp. on $S_l^+(\delta)$) belong to a small interval. Using the bounds (44)-(45) we can find these useful bounds for the difference of the Hamiltonian:

$$(-1 + s - 2s\delta)a_d|\partial_b B_{2l}| \leq H_{\partial_b B_{2l}}(\sigma_1) - H_{\partial_b B_{2l}}(\sigma_2) \leq (-1 + s + 2\delta)a_d|\partial_b B_{2l}| \tag{46}$$

for every $\sigma_1 \in S_l^-(\delta)$, $\sigma_2 \in S_l^+(\delta)$. Now we choose $s \in (0, 1)$ and $\delta \in (0, 1/2)$ such that $(-1+s-2s\delta)$ and $(-1+s+2\delta)$ will be negative. In the end, see Theorem 4.8, we will prove that the different conditions on the parameters $\delta$ and $s$ can be simultaneously satisfied. One has $|\partial_b B_{2l}| = C_d l^{d-1} = d2^d l^{d-1}$; so, using inequality (46) we find

$$e^{\beta(-1+s-2s\delta)a_d C_d l^{d-2}} \leq \frac{\pi_\beta(\Phi(\sigma))}{\pi_\beta(\sigma)} \leq e^{\beta(-1+s+2\delta)a_d C_d l^{d-2}} \quad \text{for all } \sigma \in S_l^+(\delta). \tag{47}$$

So we deduce from the first inequality of (47) that

$$\limsup_{l \to \infty} \frac{\pi_\beta(S_l^-(\delta))}{\pi_\beta(S_l^+(\delta))} = 0, \tag{48}$$

this means that for large $\beta$ -where we can apply Lemma 4.4 and Remark 4.5- we obtain:

$$\lim_{l \to \infty} \pi_\beta(S_l^+(\delta)) = 1. \tag{49}$$

Using the second inequality in (47) we deduce that

$$\limsup_{l\to\infty} \frac{\pi_\beta(S_l^+(\delta))}{\pi_\beta(S_l^-(\delta))e^{K_2\beta l^{d-1}}} = 0 \tag{50}$$

where $K_2$ can be chosen as any positive constant larger than $d2^d(1-s-2\delta)a_d = d2^{d-1}(1-s-2\delta)\rho_d\delta$, we take $K_2 = d2^d(1-s-2\delta)\rho_d\delta$. This concludes the proof. □

**Lemma 4.7.** *Let's $\delta \in (0,1/2)$. There exists $\beta_3(\delta)$ such that for all $\beta > \beta_3(\delta)$ the probability that either*

$$\sigma(t) \in S_l^+(\delta) \text{ and } \sigma(t+1) \in S_l^-(\delta) \tag{51}$$

*or*

$$\sigma(t) \in S_l^-(\delta) \text{ and } \sigma(t+1) \in S_l^+(\delta) \tag{52}$$

*is smaller than $e^{-K_3\beta l^{d-1}}$ where $K_3 = s\rho_d\delta d2^{d-5}$.*

**Proof**. The first step SW1) in the evolution of configuration of edge $\eta$ is really simple because all the variables $\{\eta_e\}_{e\in E_\Lambda}$ are independent. The probability to have on an edge $e = \{v,w\} \in E_\Lambda$ $\eta_e(t+1) = 1$ depends only on the spins $\sigma_v(t)$ and $\sigma_w(t)$. More precisely if $J_{v,w}\sigma_v(t)\sigma_w(t) > 0$ then $\eta_{v,w}$ takes value 1 with probability $1-e^{-\beta|J_{v,w}|}$; if $J_{v,w}\sigma_v(t)\sigma_w(t) \le 0$ $\eta_{v,w}$ is equal to 0 with probability 1.

Let us consider the case $\{\sigma(t) \in S_l^-(\delta)\} \to \{\sigma(t+1) \in S_l^+(\delta)\}$. Since $\sigma(t)$ is in $S_l^-(\delta)$ there are at least $(1-\delta/2)|\partial_b B_{2l}|$ edges $e = \{v,w\}$ on $\partial_b B_{2l}$ having $\sigma_v\sigma_w = -1$. Among these edges at least

$$(1-\delta/2)|\partial_b B_{2l}| - |\partial_b B_{2l}|/2$$

are negative edges. So at least $(1-\delta)|\partial_b B_{2l}|/2$ are negative edges having probability $1-p_e$ to be occupied (*i.e.* $\eta_e = 1$). There is the following uniform lower bound for the probability that one of these edges is occupied:

$$p_e \ge \inf\{1-e^{-\beta|J_e|} : e \in \partial_b B_{2l} \text{ is a negative edge }\} \ge 1-e^{-\beta s a_d}. \tag{53}$$

Now we use the standard large deviation estimation for a Bernoulli distribution; let us suppose that $\{X_i\}_{i\in\mathbb{N}}$ are i.i.d. random variable with $P_p(X_1 = 0) = 1 - p$ and $P_p(X_1 = 1) = p$ then

$$P_p\left(\sum_{i=1}^{N} X_i < xN\right) \leq e^{-I(x,p)N} \text{ if } x \leq \mathbf{E}(X_1) = p, \tag{54}$$

where the *rate function* $I(x,p) = x\log(\frac{x}{p}) + (1-x)\log(\frac{1-x}{1-p})$. For any fixed $x \in (0,1)$

$$\lim_{p\to 1^-} \frac{-I(x,p)}{(1-x)\log(1-p)} = 1. \tag{55}$$

So for any fixed $x \in (0,1)$ there is $p \in (0,1)$ large enough such that

$$I(x,p) > \frac{1-x}{2}\log\frac{1}{1-p}. \tag{56}$$

Now we can estimate the probability that in the step SW1) of SW dynamics there are less than $\frac{1}{2}\frac{1-\delta}{2}|\partial_b B_{2l}|$ occupied edges when at least $\frac{1-\delta}{2}|\partial_b B_{2l}|$ edges can be occupied

$$P_p\left(\sum_{e\in A} \eta_e < \frac{1}{2}\frac{1-\delta}{2}|\partial_b B_{2l}|\right) \leq \exp(-I(\frac{1}{2}, 1-e^{-\beta s a_d/2})\frac{1-\delta}{2}|\partial_b B_{2l}|), \tag{57}$$

we are using (54) with $p$ that is any positive number less or equal to $\inf\{1-e^{-\beta|J_e|} : e \in \partial_b B_{2l}$ is a negative edge $\}$, we have chosen $p = 1-e^{-\beta s a_d/2}$; $A \subset \partial_b B_{2l}$ is the set of edges that can be occupied and $|A| \geq \frac{1-\delta}{2}|\partial_b B_{2l}|$. Now, for a fixed $\delta \in (0,1/2)$ and for $\beta$ large enough, using (56) we can bound the probability in (57) with

$$P_p\left(\sum_{e\in A} \eta_e < \frac{1-\delta}{4}|\partial_b B_{2l}|\right) \leq \exp(-s a_d \beta|\partial_b B_{2l}|/32). \tag{58}$$

We notice that if there are at least $\frac{1-\delta}{4}|\partial_b B_{2l}|$ occupied negative edges in the step SW1) then in the second step SW2) for every occupied negative edge $e = \{v,w\}$ the product $\sigma_v\sigma_w$ will be equal to $-1$, so $K(\partial_b B_{2l}, \sigma(t))$, see definition (11), will be less than $(1+\delta)/2$ and the configuration $\sigma(t+1)$ is not in $S_l^+(\delta)$ when $\delta \in (0,1/2)$. In this case the constant $K_3$ can be chosen smaller or equal to $s a_d d 2^{d-4}$ (see inequality (58)).

The other case $\{\sigma(t) \in S_l^+(\delta)\} \to \{\sigma(t+1) \in S_l^-(\delta)\}$ is analogous and in both transitions we can take $K_3 = s a_d d 2^{d-4} = s\rho_d \delta d 2^{d-5}$ (we have used (29) for the last equality). $\square$

In the following theorem we will denote with $L$ the length of the side of the box $\Lambda = \Lambda_L$ and $l(L)$ will denote an opportune function such that one can find, for $L$ large enough, inside $\Lambda_L$ a box $B_{4l(L)}$ verifying all the hypotheses on the interactions written in the second section. It is clear that a strictly increasing function $l(L)$ exists because of finding the right interactions is an event of probability larger than zero and the interactions are *i.i.d.* random variables. In the next theorem we denote with $\Lambda_L$ the box $\mathbb{Z}^d \cap [-L, L]^d$ and $\tilde{\beta}(\delta) = \max\{\beta_1(\delta), \beta_2(\delta), \beta_3(\delta)\}$. In the next theorem we consider $\beta$ larger than $\tilde{\beta}(\delta)$ so are verified simultaneously all the inequalities of previous lemmas.

**Theorem 4.8.** *Let $\delta \in (0, 1/2)$, $\beta > \tilde{\beta}(\delta)$ , and $f_{L,\delta} : \Omega_{\Lambda_L} \to \mathbb{R}$*

$$f_{L,\delta}(\omega) = \mathbf{1}_{S^-_{l(L)}(\delta)}(\omega) \tag{59}$$

*then, the spectral radius $R_L$ goes to one when $L$ increase to infinity $Q_J$-almost surely and so*

$$\lim_{L\to\infty} \tau_{exp,L} = \infty \quad Q_J - a.s. \tag{60}$$

**Proof**. First we need to show that it is possible to select the constant $K_1$, $K_2$ and $K_3$ of Lemmas 4.4-4.7 and Remark 4.5 such that they satisfy the inequalities $0 < K_2 < K_1$ and $K_2 < K_3$. We show that there are some values of the parameters $\delta$ and $s$ such that all the inequalities between $K_1$, $K_2$ and $K_3$ are verified. The inequalities to be verified are:

$$\begin{cases} 1 - s - 2\delta > 0 \qquad (\Rightarrow K_1 > 0) \\ \delta^2\rho_d^2 > 16(1 - s - 2\delta) \qquad (\Rightarrow K_1 > K_2) \\ s\rho_d\delta > 64(1 - s - 2\delta) \qquad (\Rightarrow K_3 > K_2) \\ 0 < s < 1 \\ 0 < \delta < 1/2 \end{cases} \tag{61}$$

We remind that $\rho_d$ is a fixed geometrical constant in the open interval $(0, 1)$ (see Lemma 4.2 and Lemma 4.3); it is immediate to see that for every $\delta \in (0, 1/2)$ there are some values of the parameter $s$ such that all the inequalities are verified (a non empty interval of values). One can choose any value of $s$ in the open interval

$$\left(\max\left\{\frac{64}{\rho_d\delta + 64}(1 - 2\delta), 1 - 2\delta - \frac{\delta^2\rho_d}{16}\right\}, 1 - 2\delta\right);$$

$s < 1 - 2\delta$ from the first inequality in (61), $s > \max\left\{\frac{64}{\rho_d\delta+64}(1-2\delta), 1-2\delta-\frac{\delta^2\rho_d}{16}\right\}$ from the second and third inequalities in (61). It is immediate to see that also the other inequalities are verified and

$$\max\left\{\frac{64}{\rho_d\delta+64}(1-2\delta), 1-2\delta-\frac{\delta^2\rho_d}{16}\right\} < 1-2\delta$$

for every $\delta \in (0, 1/2)$. In the rest of the proof we will omit the dependence of the eigenvalues from the parameter $\delta$.

Now we use the Dirichlet's formula (see equation (5)) and the variational characterization of $\lambda_1$ to prove that the spectral gap goes to zero when $L$ increase to infinity. The usual variational characterization of the second eigenvalue gives:

$$\lambda_1 = 1 - \inf\left\{\frac{C_{ff}(1)}{Var_\pi(f)} : 0 < Var_\pi(f) < \infty\right\} \tag{62}$$

and also

$$\lambda_1^t = 1 - \inf\left\{\frac{C_{ff}(t)}{Var_\pi(f)} : 0 < Var_\pi(f) < \infty\right\} \tag{63}$$

for every $t \in \mathbb{N}$; where $t$ is the number in which the configuration is updated with SW algorithm.

If we can prove that there exists an increasing sequence of times $t_L = t(l(L)) \uparrow \infty$ such that

$$\lim_{L\to\infty} \frac{C_{f_L,f_L}(t_L)}{Var_\pi(f_L)} = 0 \quad Q_J - a.s. \tag{64}$$

then we get that the spectral gap goes to zero $Q_J$-a.s. (see for similar arguments [9]).

Let us define $\tilde{f}_l = f_L = \mathbf{1}_{S_l^-(L)}$; then for every time $t \in \mathbb{N}$ we have the following upper bound

$$C_{f,f}(t) \leq t\pi_\beta(\Omega \setminus (S_l^- \cup S_l^+)) + t \max_{\omega\in S_l^-,\omega'\in S_l^+} M_\beta^{SW}(\omega,\omega') + t \max_{\omega\in S_l^+,\omega'\in S_l^-} M_\beta^{SW}(\omega,\omega') \tag{65}$$

where we are supposing that there is a box $B_{4l}$ with all the interactions described in the second section. The first summand in RHS of (65) is an upper bound for the probability that in the time $\{1,\dots,t\}$ the chain visits a configuration out of $S_l^+ \cup S_l^-$; this is a consequence of the initial measure that is equal to the stationary measure $\pi_\beta$, so at any time the associated measure of the Markov chain is equal to $\pi_\beta$; the upper bound is a consequence of sub-additivity. The second summand in RHS of (65) is an upper bound of a transition from a state of $S_l^-$ to a state of $S_l^+$. Analogously the third summand in (65) control the probability of a transition from $S_l^+$ to $S_l^-$. We are imposing that every time that there is a transition in $s \in \{1, \ldots, t\}$ from $S_l^-$ to $S_l^+$ or a transition from $S_l^+$ to $S_l^-$ or the chain visits a configuration in $\Omega \setminus (S_l^+ \cup S_l^-)$ then the function $(\mathbf{1}(S_l^-)(\omega_t) - \mathbf{1}(S_l^-)(\omega_0))^2$ of formula (5) is equal to 1 where $\omega_s$ is the configuration of the process at time $s \in \{1, \ldots, t\}$ (this is an obvious upper bound). To find an upper bound for the RHS of (65) we will use Lemma 4.4 and Lemma 4.7, so one obtains that

$$C_{f_l,f_l}(t) \leq t[e^{-K_1\beta l^{d-1}} + 2e^{-K_3\beta l^{d-1}}] \tag{66}$$

For $C_{f_l,f_l}(0) = Var(f)$, using Lemma 4.6, we obtain the following lower bound

$$C_{f_l,f_l}(0) \geq \frac{1}{2}e^{K_2\beta l^{d-1}}.$$

So

$$\frac{C_{f_l,f_l}(t)}{C_{f_l,f_l}(0)} \leq 2t[e^{-(K_1-K_2)\beta l^{d-1}} + 2e^{-(K_3-K_2)\beta l^{d-1}}], \tag{67}$$

where the constant $(K_1 - K_2)$ and $(K_3 - K_2)$ are positive. So one can choose an increasing sequence of times $t_l$, depending on $l$, such that $t_l \rightarrow_{l\rightarrow\infty} \infty$ and

$$\lim_{l\rightarrow\infty} \frac{C_{f_l,f_l}(t_l)}{C_{f_l,f_l}(0)} = 0.$$

To finish the proof we stress that one can find a sequence of boxes $\Lambda_L$ (with side length $L$) such that just a finite number of these boxes $\Lambda_L$ have not inside a square $B_{4l}$ with all the interactions verifying hypotheses 1-3 of the second section. So $R_L$ goes to one $Q_J$-a.s. and the exponential autocorrelation time goes to infinity using the relation (11). □

Some remarks: We can consider a limit dynamics when the temperature parameter $T$ goes to zero; it is clear that in this case there is a different behavior for the ferromagnetic systems with respect to the frustrated ones. For the ferromagnetic systems all the configurations that are not ground states for the Hamiltonian are transient states and the unique invariant measure of the associated Markov chain is uniform on both the ground states (all the spins equal to 1 or -1). Moreover the ground states form an irreducible class so with probability one the Markov chain visits all these ground states (the states are finite because we are considering just a finite box $\Lambda$). For frustrated systems there are different invariant measures and the ground states are divided in several irreducible classes; so it is immediate to see that in this case the spectral gap is equal to zero because there are states that are not reachable, starting from some states, in spite of a positive probability for some invariant measures.

We want to stress that not only our proof does not work for the ferromagnetic systems but also we believe that really there is a different behavior for all the small values of the temperature between ferromagnetic and frustrated systems. Our conjecture is that for every ferromagnetic system there exists $T_0$ such that for each $T \in [0, T_0)$ the spectral gap of the associated Markov chain is uniformly bounded away from zero independently of the size of the systems (the size of the box $\Lambda$). This is not in contrast with [18] because a Markov process on an infinite space of states can be non-ergodic also if all the terms in the sequence of its approximations on finite spaces of states have uniform bounds (away from zero) for the spectral gap.

Vice versa in this paper we have just proved that for every frustrated disordered system, where the interactions have $0 \in supp \ \ Q_J$, the spectral gap shrinks to zero when the size of the system increases to infinity and the temperature is small enough. So we believe that there is a different behavior for small values of the temperature between the ferromagnetic and the disordered not ferromagnetic systems.

Following [9] one could prove that for small and large values of $\beta$ there is a different behavior of the convergence to the equilibrium with respect to the distance $D_\rho$. The distance $D_\rho$ between two measures $\mu_1$ and $\mu_2$ defined on $\{-1,1\}^\Lambda$ is:

$$D_\rho(\mu_1, \mu_2) = \sup_{i \in \Lambda} \sum_{A, i \in A \subset \Lambda} \sum_{\sigma_A \in \{-1,1\}^A} |\mu_1(\sigma_{A+i}) - \mu_2(\sigma_{A+i})| \, \rho^{|A|} \tag{68}$$

where $A$ is a connected vertex set; $\rho \in (0,1)$ is a parameter. For small values of $\beta$ and of the parameter $\rho$ we obtain that

$$D_\rho(\mu_0 M_\Lambda^{SW^t}, \mu_{\Lambda,J,\beta}) < M e^{-ct} \tag{69}$$

where $c > 0$ can be chosen independent of $\Lambda$. For large values of $\beta$ the relation (69) is not true for any fixed value of $\rho$.

## 5. Conclusions

We have proved some bounds on the spectral radius of the Swendsen-Wang dynamics and we found a different behavior for large or small values of the temperature parameter. Our upper bound on the spectral gap is not explicit but if one needs an explicit evaluation one should only perform our calculation in Theorem 4.8 preserving all the coefficients. Also one could find, for every distribution of the interactions $\mathcal{J}$, the optimal values of the parameters $\delta$ and $s$ to estimate, using the Dirichlet's formula, the spectral gap of $M^{SW}_{\Lambda,\beta}$.

A comparison with the SW dynamics of the ferromagnetic Ising model would be interesting. We believe that the ferromagnetic systems carry out better performances than the disordered not ferromagnetic one at every value of the temperature $1/\beta$; moreover we think that the spectral gap in the ferromagnetic case does not shrink to zero when the size increases to infinity if the temperature is small enough.

**Acknowledgments**. *I wish to thank Charles Newman and Alan Sokal for beneficial conversations and for encouraging me to work on this paper. I thank Jeff Cheeger which have given me some ideas to prove Lemma 4.2.*

Università di Roma *La Sapienza* Dipartimento di Matematica “Guido Castelnuovo”

Piazzale Aldo Moro 2, 00185 Roma, Italia